%% file: main.tex
\documentclass[letterpaper]{article}

\input{packages}
\input{environment}

\title{Mixed-variable policy-based optimization}

\def\size{7.3cm}
\author{
  \parbox{\size}{\centering J. Viquerat\thanks{Corresponding author}}\\
  MINES Paristech, CEMEF\\
  PSL - Research University\\
  \texttt{jonathan.viquerat@minesparis.psl.eu}
}

\begin{document}
\newgeometry{left=3.5cm,right=3.5cm,top=2.5cm,bottom=3cm}
\maketitle

\begin{abstract}
The optimization of mixed-variable problems remains a significant challenge. We propose an extension of the policy-based optimization method that handles mixed-variables problems in a natural way, through a simple policy combination. This is achieved by independently sampling from a multivariate normal distribution for the continuous domain, and from multiple categorical distributions for the discrete choices. Results demonstrate that the agent successfully yields high-quality solutions on a classical problem of electromagnetics, showcasing its robustness.
\end{abstract}


\section{Introduction}

Mixed-variable optimization addresses problems containing both continuous and discrete decision variables. The inherent difficulty arises from a "worst of both worlds" structure: a combinatorial, non-convex search space induced by the discrete variables, coupled with the often complex, non-linear landscape of the continuous variables. This hybrid nature invalidates many classical optimization techniques, necessitating specialized algorithms. The standard methods for tackling these problems can be broadly categorized into two families, distinguished by the information available about the objective function.

For white-box problems (where a precise mathematical formulation is known), the outer approximation approach is a fundational, well-established method \cite{Fletcher1994}. This method iterates between solving continuous nonlinear subproblems and a simplified linear master problem that guides the selection of discrete variables \cite{Duran1986}. For black-box problems, where only function evaluations are possible, a diverse range of heuristic and derivative-free methods are employed. Evolutionary algorithms are a natural fit, as they can handle mixed-variable solutions through specialized genetic representations and operators tailored to each variable type \cite{Coello2002}. Other population-based metaheuristics, such as particle swarm optimization, have also been successfully adapted by defining specific update rules or mappings for the discrete components of a solution \cite{Laskari2002}. More recently, bayesian optimization has emerged as a powerful, model-based approach for expensive black-box problems. This method relies on building a statistical surrogate model of the objective function, using specialized kernels capable of handling the hybrid continuous-discrete space to intelligently guide the search \cite{Garrido2020}.

Here, we propose an extension of the policy-based optimization (PBO) method, introduced in \cite{Viquerat2023}, which naturally handles the mixed nature of the problem through a combination of policies. PBO is a black-box optimization algorithm based on a single-step deep reinforcement learning (DRL) approach. The method uses independent neural networks to learn the parameters of a multivariate normal search distribution. It operates on the principle of single-step episodes, where the policy is independent of the state and aims to maximize an instantaneous reward. In essence, the algorithm acts as an evolution strategy where the update rules of the distribution are learned through the policy gradient-based training of neural networks instead of being defined by analytical heuristics.

In the following, we first introduce a revisited version of PBO, with core modifications to the original algorithm. Then, the mixed-variable PBO approach is described. Finally, the method is used to optimize a stacked dielectric reflector, a sandwich-like structure of dielectric materials aiming a providing maximal reflectance over a range of wavelengths. In this context, the method is successfully tested up to 40 degrees of freedom.

\section{Revisiting PBO}

In this section, we revisit the PBO method introduced in \cite{Viquerat2023} with two major modifications for improved stability and performance.

\paragraph{Mean update} We noticed that the use of a neural network to generate the mean vector $\V{m}$ could sometimes lead to unstable updates due to its reduced update buffer size. In the present version, it is replaced with a weighted recombination of elite points from the previous generation:

\begin{equation*}
  \V{m}_{t+1} = \sum_i w^i x^i_t, \text{ with } w_i = \frac{\hat{R}_i}{\sum_j \hat{R}_j},
\end{equation*}

where the sums are performed on the selected elite points from last generation, and $\hat{R}_i$ is the $i^\text{th}$ component of the whitened reward vector, as described in \cite{Viquerat2023}. The update and generation of the covariance matrix remains unchanged.

\paragraph{Elite points} At each generation, the original PBO approach pruned negative advantage samples, leading to possibly unbalanced update buffers. In this updated method, we retain a fixed amount of elite points at each generation, usually equal to half the number of sampled points.

\section{Mixed-variable PBO}

The introduction of discrete variables in the original continuous PBO formulation is a straightforward task. From this point, we assume the considered problem holds $n_c$ continuous variables, and $n_d$ discrete variables. For each discrete variable, we denote as $d_k$ the number of categories this variable can take. An additional network is added to represent the discrete variables, with a final linear layer of size $\sum_k d_k$, outputting the logits of all possible categories for all discrete variables at once, which are then used to parameterize the related categorical distributions. At sampling time, the output vector of size $\sum_k d_k$ is split into $n_d$ vectors of adequate sizes, and a categorical distribution is generated for each discrete variable, from which samples are then drawn. This defines a joint policy $\pi_{\theta_d}$ for all discrete variables:

\begin{equation*}
  \pi_{\theta_d} (a_d) = \prod_{i=1}^{n_d} \pi_{\theta_d} (a_{d,i}).
\end{equation*}

Given the independence of the continuous $\pi_{\theta_c}$ and discrete $\pi_{\theta_d}$ policies, the log-probability of a complete mixed-variable action $a=(a_c, a_d)$ is the sum of the individual log-probabilities:

\begin{equation*}
  \log \pi_{\theta} (a) = \log \pi_{\theta_c} (a) + \log \pi_{\theta_d} (a),
\end{equation*}

where $\theta$ represents the full set of parameters for all policy networks. Thanks to this simple combination, the computation of the advantage vector and the loss function remain similar to what is described in \cite{Viquerat2023}.

\section{Results}

We consider the optimization of a multi-layer dielectric mirror, a fundamental challenge in optical engineering. The goal is to select a sequence of thin, transparent materials and their respective thicknesses, in order to maximize the reflection of light over a specified range of wavelengths. The main difficulty lies in the complex physics of wave interference; as light passes through the stack, reflections from each interface interfere with one another in ways that are highly sensitive to the thickness and refractive index of each layer. This creates a high-dimensional and non-convex optimization landscape with multiple local optima, making it difficult to find a globally optimal design that is both highly reflective and free of unwanted reflection dips.

The computation of the reflectance is performed using the transfer matrix method. The transfer matrix method is a powerful analytical tool used in optics to calculate the transmission and reflection of light through a stack of thin-film layers. The method characterizes the effect of each layer and interface on the wave's electric and magnetic fields using a specific $2 \times 2$ transfer matrix. From the resulting total system matrix, which is found by multiplying the individual matrices in sequence, macroscopic properties such as the overall reflectance and transmittance of the stack can be directly and efficiently calculated \cite{Hecht2017}.

For our specific implementation, the task is to design a stack of $n_l$ layers, choosing from two materials: titanium dioxyde ($\text{TiO}_2$), with a refractive index $n_{\text{TiO}_2} = 2.4$, and magnesium fluoride ($\text{MgF}_2$), of refractive index $n_{\text{MgF}_2} = 1.38$. These layers are deposited on a glass substrate ($n = 1.52$) and are surrounded by air ($n = 1.0$). The optimization algorithm controls two sets of variables: a discrete choice of material for each layer, and a continuous thickness for each layer, typically bounded between $50$ nm and $150$ nm. The performance of a given design is evaluated by calculating its average reflectance across a target spectrum, a band between $300$ nm and $500$ nm. In the following, we consider $n_l = 20$, implying the optimization problem has $40$ degrees of freedom.

\paragraph{Maximizing reflectance} First, we evaluate the method the simplest objective, \textit{i.e.} maximizing the average reflectance over the considered bandwidth. The cost function we minimize is therefore:

\begin{equation*}
  c(\V{x}) = - \left< \rho (\lambda, \V{x}) \right>_{\lambda},
\end{equation*}

where $\V{x}$ is the design vector. The results are presented in figure \ref{fig:max}. As can be seen, the proposed design creates a reflectance spectrum with many flat regions, with an average reflectance of $0.9306$ over the considered spectrum. Yet, this is achieved  at the cost of several reflectance dips in the spectrum. Convergence is reached in approximately $7.5k$ evaluations, with a moderate standard deviation over the different runs.

\input{fig/max}

\paragraph{Maximally-flat reflectance} We now consider a modified cost function, in which we include a penalization term for the performance dips appearing in figure \ref{fig:max}. The objective is to find a trade-off between the flatness of the reflectance over the considered spectrum and its mean value:

\begin{equation*}
  c(\V{x}) = - \left< \rho (\lambda, \V{x}) \right>_{\lambda} + \alpha (\max_{\lambda} \rho(\lambda, \V{x}) - \min_{\lambda} \rho(\lambda, \V{x})).
\end{equation*}

In the following, we choose $\alpha = 0.1$, and the results are presented in figure \ref{fig:maxflat}. As can be observed, the difference between the maximal and minimal reflectance values is reduced compared to the previous case, at the expanse of (i) a reduced mean reflctance of $0.907$, (ii) a reduced flatness of the spectrum, and (iii) an increased variance of the cost function over the different runs, indicating a more complex cost function landscape for the optimizer, probably linked to the existence of multiple local minima.

\input{fig/flat}

\section{Conclusion}

In this work, a straightforward yet effective extension of the policy-based optimization method is presented to address mixed-variable optimization problems. The inherent difficulty of these problems stems from their hybrid nature, which combines a combinatorial discrete search space with a complex continuous landscape. The proposed framework naturally navigates this challenge by integrating separate policies for the continuous and discrete components of a solution. This is achieved by sampling continuous variables from a multivariate normal distribution while concurrently drawing discrete variables from categorical distributions. The effectiveness of the method was evaluated on the design of a 20-layer dielectric mirror with 40 degrees of freedom.

\bibliographystyle{unsrt}
\bibliography{bib}

\end{document}

%% file: packages.tex
\usepackage[utf8]{inputenc}
\usepackage{arxiv}
\usepackage{float}
\usepackage{enumitem}
\usepackage[english]{babel}
\usepackage{amssymb}
\usepackage{amsmath,mathtools}
\usepackage{stmaryrd} 								
\usepackage[table]{xcolor}
\usepackage{graphicx}
\usepackage{algorithmicx}
\usepackage{algpseudocode}
\usepackage{algorithm}
\usepackage{multirow}
\usepackage{cprotect}
\usepackage{bbold}
\usepackage{overpic}
\usepackage{psfrag}
\usepackage{bm}
\usepackage{lineno}
\usepackage{url}
\usepackage{physics}								
\usepackage{epsfig}
\usepackage{marvosym}
\usepackage{verbatim}
\usepackage{subcaption}
\usepackage{caption}
\usepackage{hyperref}
\hypersetup{colorlinks=true,linkcolor=red,citecolor=blue}
\usepackage{ulem}
\usepackage{lmodern}
\usepackage{booktabs} 
\usepackage{moresize} 
\usepackage{siunitx}
\usepackage{mathtools}
\usepackage{tikz}
\usepackage{pgfplots}
\usepackage{pgfplotstable}
\pgfplotsset{compat=1.14}
\usepgfplotslibrary{fillbetween}
\usetikzlibrary{shapes, arrows, arrows.meta, calc, positioning, decorations,
decorations.pathmorphing, decorations.text, spy}
\pgfplotsset{translate gnuplot=true}
\usepackage{etex} 								
\usepgfplotslibrary{external} 						
\tikzsetexternalprefix{ext/}							
\usepackage{booktabs}
\usepackage{makecell}
\usepackage{rotating,tabularx}

%% file: environment.tex
\definecolor{blue1}		{RGB}{0,177,234}				
\definecolor{blue2}		{RGB}{76,200,239}				
\definecolor{blue3}		{RGB}{127,215,244}				
\definecolor{blue4}		{RGB}{178,231,248}				
\definecolor{bluegray1}	{RGB}{0,127,167}				
\definecolor{bluegray2}	{RGB}{76,165,193}				
\definecolor{bluegray3}	{RGB}{127,191,211}				
\definecolor{bluegray4}	{RGB}{178,216,228}				
\definecolor{gray1}		{RGB}{76,84,93}				
\definecolor{gray2}		{RGB}{129,135,141}				
\definecolor{gray3}		{RGB}{165,169,174}				
\definecolor{gray4}		{RGB}{201,203,206}				
\definecolor{gray5}		{RGB}{230,230,230}				
\definecolor{gray6}		{RGB}{245,245,245}				
\definecolor{orange1}	{RGB}{255,126,46}				
\definecolor{orange2}	{RGB}{255,164,108}				
\definecolor{orange3}	{RGB}{255,190,150}				
\definecolor{orange4}	{RGB}{255,216,192}				
\definecolor{brown1}		{RGB}{205,133,63}				
\definecolor{green1}		{RGB}{0,168,107}				
\definecolor{green2}		{RGB}{30,198,137}				
\definecolor{green3}		{RGB}{60,228,167}				
\definecolor{green4}		{RGB}{90,255,197}				
\definecolor{purple1}		{RGB}{89,89,171}				
\definecolor{purple2}		{RGB}{159,159,195}				
\definecolor{purple3}		{RGB}{189,189,225}				
\definecolor{purple4}		{RGB}{209,209,255}				
\definecolor{teal1}		{RGB}{0,128,128}				
\definecolor{teal2}		{RGB}{100,168,168}				
\definecolor{teal3}		{RGB}{130,198,198}				
\definecolor{teal4}		{RGB}{160,228,228}				
\definecolor{red1}		{RGB}{255,0,0}					
\definecolor{red2}		{RGB}{255,30,30}				
\definecolor{red3}		{RGB}{255,60,60}				
\definecolor{red4}		{RGB}{255,90,90}				
\definecolor{royal1}		{RGB}{2,119,189}				
\definecolor{royal2}		{RGB}{32,149,219}				
\definecolor{royal3}		{RGB}{62,179,249}				
\definecolor{royal4}		{RGB}{92,209,255}				

\pgfplotsset{
    colormap={custom_map}{[5pt]
            rgb255(0pt)=(255,126,46);
            rgb255(500pt)=(255,190,150);
            rgb255(1000pt)=(0,177,234);
            rgb255(1500pt)=(127,215,244);
    },
}

\renewcommand{\emph}[1]{\textit{#1}}
\setenumerate{label=\textcolor{bluegray1}{$\bm{\diamond}$},itemsep=0.0pt,topsep=5pt}

\newcommand{\V}[1]{\bm{#1}}

%% file: fig/max.tex
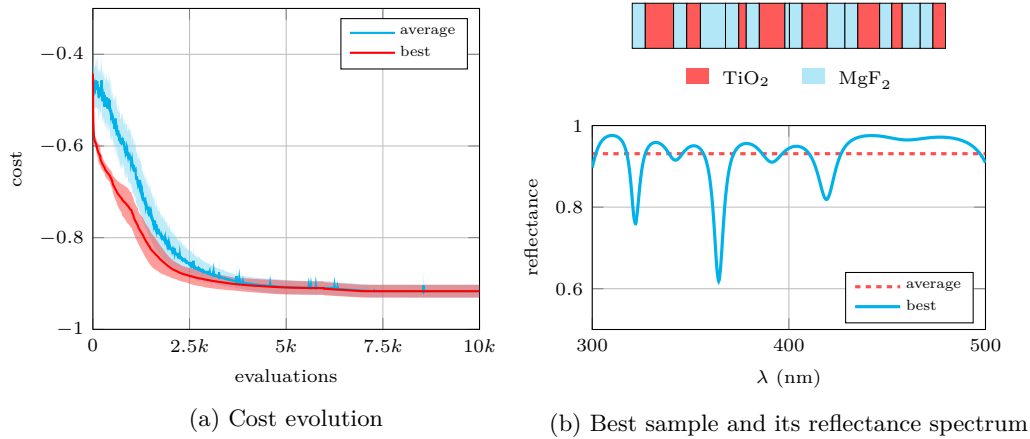
\begin{figure}
\centering
\begin{minipage}{.45\linewidth}
\begin{subfigure}[t]{.99\textwidth}
\centering
\begin{tikzpicture}[
      trim axis left, trim axis right, font=\scriptsize,
      upper/.style={	name path=upper, smooth, draw=none},
      lower/.style={	name path=lower, smooth, draw=none},]
\begin{axis}[	xmin=0, xmax=10000, scale=0.75,
              ymax=-0.3, ymin=-1.0,
              scaled x ticks=false,
              xtick={0, 2500, 5000, 7500, 10000},
              xticklabels={$0$,$2.5k$,$5k$,$7.5k$,$10k$},
              legend cell align=left, legend pos=north east,
              legend style={nodes={scale=0.8, transform shape}},
              every tick label/.append style={font=\scriptsize},
              grid=major, xlabel=evaluations, ylabel=cost]

      \legend{average, best}

      \addplot[upper, forget plot]                   table[x index=0,y index=3] {fig/max_avg.dat};
      \addplot[lower, forget plot]                   table[x index=0,y index=2] {fig/max_avg.dat};
      \addplot[fill=blue3, opacity=0.5, forget plot] fill between[of=upper and lower];
      \addplot[draw=blue1, thick, smooth]            table[x index=0,y index=1] {fig/max_avg.dat};

      \addplot[upper, forget plot]                  table[x index=0,y index=6] {fig/max_avg.dat};
      \addplot[lower, forget plot]                  table[x index=0,y index=5] {fig/max_avg.dat};
      \addplot[fill=red3, opacity=0.5, forget plot] fill between[of=upper and lower];
      \addplot[draw=red1, thick, smooth]		        table[x index=0,y index=4] {fig/max_avg.dat};

\end{axis}
\end{tikzpicture}
\caption{Cost evolution}
\label{fig:max_cost}
\end{subfigure}
\end{minipage}
\begin{minipage}{.45\linewidth}
\begin{subfigure}[t]{.99\textwidth}
\centering
\begin{tikzpicture}
\pgfplotstableread[]{
thickness matindex
69.2685 1
150.289 0
69.8393 1
72.822 0
133.232 1
70.3743 1
39.8179 0
69.4554 1
135.72 0
24.1919 1
67.4575 1
133.263 0
94.0047 1
68.5537 1
116.029 0
64.3018 1
55.6979 0
95.5128 1
67.4992 1
66.764  0
}\stackdata

  \pgfmathsetmacro{\x}{0}
  \pgfmathsetmacro{\scaling}{0.0025}
  \pgfplotstablegetrowsof{\stackdata}
  \pgfmathtruncatemacro{\numrows}{\pgfplotsretval}

  \foreach \i in {0,...,\numexpr\numrows-1\relax} {
    \pgfplotstablegetelem{\i}{thickness}\of{\stackdata} \let\thickness\pgfplotsretval
    \pgfplotstablegetelem{\i}{matindex}\of{\stackdata} \pgfmathtruncatemacro{\matindex}{\pgfplotsretval}

    \ifcase\matindex
      \def\layercolor{red4}
    \or
      \def\layercolor{blue4}
    \fi

    \pgfmathsetmacro{\scaledthickness}{\thickness * \scaling}
    \draw [fill=\layercolor, draw=black, line width=0.2pt] (\x, 0) rectangle ++(\scaledthickness, 0.6);
    \pgfmathsetmacro{\newx}{\x + \scaledthickness}
    \global\let\x\newx
  }

\node[inner sep=0pt, outer sep=0pt, anchor=north] at (\x/2, -0.2cm) {
  \begin{tabular}{llll}
    \tikz\fill[red4] (0,0) rectangle (0.3cm,0.2cm); & \hspace{-10pt} \scriptsize $\text{TiO}_2$  & \tikz\fill[blue4] (0,0) rectangle (0.3cm,0.2cm); & \hspace{-10pt} \scriptsize $\text{MgF}_2$
  \end{tabular}
};
\end{tikzpicture}
\end{subfigure}

\smallskip
\smallskip

\begin{subfigure}[b]{.99\textwidth}
\centering
\begin{tikzpicture}[
      trim axis left, trim axis right, font=\scriptsize,
      upper/.style={	name path=upper, smooth, draw=none},
      lower/.style={	name path=lower, smooth, draw=none},]
\begin{axis}[	xmin=300, xmax=500, scale=0.75,
              ymin=0.5, ymax=1,
              width=\linewidth/0.76,
              height=0.8\textwidth,
              scaled x ticks=false,
              xtick={300, 400, 500},
              legend cell align=left, legend pos=south east,
              legend style={nodes={scale=0.8, transform shape}},
              every tick label/.append style={font=\scriptsize},
              grid=major, xlabel=$\lambda$ (nm), ylabel=reflectance]

      \legend{average, best}

			\addplot[very thick, opacity=0.7, dash pattern=on 2pt, draw=red1]	coordinates {(300,0.931) (500,0.931)};
      \addplot[draw=blue1, very thick, smooth] table[x index=0,y index=1] {fig/max.dat};

\end{axis}
\end{tikzpicture}
\caption{Best sample and its reflectance spectrum}
\label{fig:max_spectrum}
\end{subfigure}
\end{minipage}

\caption{\textbf{Results for the maximal reflectance problem.} (Left) The curves represent the evolution of the average and best cost so far, respectively. The optimization process is performed 5 times: the averaged version of each curve is represented by the thick line, while the shaded area represent the standard deviation. (Top right) The best stack found by the optimization process. (Bottom right) The associated reflectance spectrum.}
\label{fig:max}
\end{figure}

%% file: fig/flat.tex
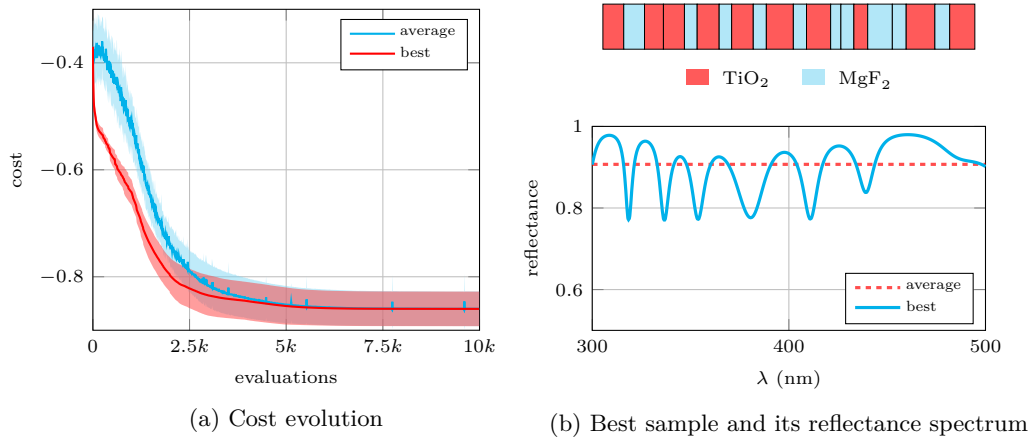
\begin{figure}
\centering
\begin{minipage}{.45\linewidth}
\begin{subfigure}[t]{.99\textwidth}
\centering
\begin{tikzpicture}[
      trim axis left, trim axis right, font=\scriptsize,
      upper/.style={	name path=upper, smooth, draw=none},
      lower/.style={	name path=lower, smooth, draw=none},]
\begin{axis}[	xmin=0, xmax=10000, scale=0.75,
              ymax=-0.3, ymin=-0.9,
              scaled x ticks=false,
              xtick={0, 2500, 5000, 7500, 10000},
              xticklabels={$0$,$2.5k$,$5k$,$7.5k$,$10k$},
              legend cell align=left, legend pos=north east,
              legend style={nodes={scale=0.8, transform shape}},
              every tick label/.append style={font=\scriptsize},
              grid=major, xlabel=evaluations, ylabel=cost]

      \legend{average, best}

      \addplot[upper, forget plot]                   table[x index=0,y index=3] {fig/flat_avg.dat};
      \addplot[lower, forget plot]                   table[x index=0,y index=2] {fig/flat_avg.dat};
      \addplot[fill=blue3, opacity=0.5, forget plot] fill between[of=upper and lower];
      \addplot[draw=blue1, thick, smooth]            table[x index=0,y index=1] {fig/flat_avg.dat};

      \addplot[upper, forget plot]                  table[x index=0,y index=6] {fig/flat_avg.dat};
      \addplot[lower, forget plot]                  table[x index=0,y index=5] {fig/flat_avg.dat};
      \addplot[fill=red3, opacity=0.5, forget plot] fill between[of=upper and lower];
      \addplot[draw=red1, thick, smooth]		        table[x index=0,y index=4] {fig/flat_avg.dat};

\end{axis}
\end{tikzpicture}
\caption{Cost evolution}
\label{fig:maxflat_cost}
\end{subfigure}
\end{minipage}
\begin{minipage}{.45\linewidth}
\begin{subfigure}[t]{.99\textwidth}
\centering
\begin{tikzpicture}
\pgfplotstableread[]{
thickness matindex
110.746 0
110.376 1
100.482 0
111.693 0
67.618  1
116.842 0
66.2679 1
114.859 0
69.4919 1
143.449 0
69.8565 1
128.275 0
54.3696 1
68.9895 1
73.5527 0
130.846 1
73.6165 1
152.288 0
78.8352 1
133.578 0
}\stackdata

  \pgfmathsetmacro{\x}{0}
  \pgfmathsetmacro{\scaling}{0.0025}
  \pgfplotstablegetrowsof{\stackdata}
  \pgfmathtruncatemacro{\numrows}{\pgfplotsretval}

  \foreach \i in {0,...,\numexpr\numrows-1\relax} {
    \pgfplotstablegetelem{\i}{thickness}\of{\stackdata} \let\thickness\pgfplotsretval
    \pgfplotstablegetelem{\i}{matindex}\of{\stackdata} \pgfmathtruncatemacro{\matindex}{\pgfplotsretval}

    \ifcase\matindex
      \def\layercolor{red4}
    \or
      \def\layercolor{blue4}
    \fi

    \pgfmathsetmacro{\scaledthickness}{\thickness * \scaling}
    \draw [fill=\layercolor, draw=black, line width=0.2pt] (\x, 0) rectangle ++(\scaledthickness, 0.6);
    \pgfmathsetmacro{\newx}{\x + \scaledthickness}
    \global\let\x\newx
  }

\node[inner sep=0pt, outer sep=0pt, anchor=north] at (\x/2, -0.2cm) {
  \begin{tabular}{llll}
    \tikz\fill[red4] (0,0) rectangle (0.3cm,0.2cm); & \hspace{-10pt} \scriptsize $\text{TiO}_2$  & \tikz\fill[blue4] (0,0) rectangle (0.3cm,0.2cm); & \hspace{-10pt} \scriptsize $\text{MgF}_2$
  \end{tabular}
};
\end{tikzpicture}
\end{subfigure}

\smallskip
\smallskip

\begin{subfigure}[b]{.99\textwidth}
\centering
\begin{tikzpicture}[
      trim axis left, trim axis right, font=\scriptsize,
      upper/.style={	name path=upper, smooth, draw=none},
      lower/.style={	name path=lower, smooth, draw=none},]
\begin{axis}[	xmin=300, xmax=500, scale=0.75,
              ymin=0.5, ymax=1,
              width=\linewidth/0.76,
              height=0.8\textwidth,
              scaled x ticks=false,
              xtick={300, 400, 500},
              legend cell align=left, legend pos=south east,
              legend style={nodes={scale=0.8, transform shape}},
              every tick label/.append style={font=\scriptsize},
              grid=major, xlabel=$\lambda$ (nm), ylabel=reflectance]

      \legend{average, best}

			\addplot[very thick, opacity=0.7, dash pattern=on 2pt, draw=red1]	coordinates {(300,0.907) (500,0.907)};
      \addplot[draw=blue1, very thick, smooth] table[x index=0,y index=1] {fig/flat.dat};

\end{axis}
\end{tikzpicture}
\caption{Best sample and its reflectance spectrum}
\label{fig:maxflat_spectrum}
\end{subfigure}
\end{minipage}

\caption{\textbf{Results for the maximally flat reflectance problem.} (Left) The curves represent the evolution of the average and best cost so far, respectively. The optimization process is performed 5 times: the averaged version of each curve is represented by the thick line, while the shaded area represent the standard deviation. (Top right) The best stack found by the optimization process. (Bottom right) The associated reflectance spectrum.}
\label{fig:maxflat}
\end{figure}

%% file: main.bbl
\begin{thebibliography}{1}

\bibitem{Fletcher1994}
R.~Fletcher and S.~Leyffer.
\newblock Solving mixed integer nonlinear programs by outer approximation.
\newblock {\em Mathematical Programming}, 66(1-3):327--349, 1994.

\bibitem{Duran1986}
M.~A. Duran and I.~E. Grossmann.
\newblock An outer-approximation algorithm for a class of mixed-integer
  nonlinear programs.
\newblock {\em Mathematical Programming}, 36(3):307--339, 1986.

\bibitem{Coello2002}
C.~A. Coello~Coello.
\newblock A comprehensive survey of evolutionary-based multiobjective
  optimization techniques.
\newblock {\em Knowledge and Information Systems}, 1(3):269--308, 2002.

\bibitem{Laskari2002}
E.~C. Laskari, K.~E. Parsopoulos, and M.~N. Vrahatis.
\newblock Particle swarm optimization for integer programming.
\newblock In {\em Proceedings of the 2002 IEEE Congress on Evolutionary
  Computation (CEC 2002)}, volume~2, pages 1582--1587, 2002.

\bibitem{Garrido2020}
E.~C. Garrido-Merchán and D.~Hernández-Lobato.
\newblock Dealing with categorical and integer variables in bayesian
  optimization with gaussian processes.
\newblock {\em Neurocomputing}, 380:20--35, 2020.

\bibitem{Viquerat2023}
J.~Viquerat, R.~Duvigneau, P.~Meliga, A.~Kuhnle, and E.~Hachem.
\newblock Policy-based optimization: single-step policy gradient method seen as
  an evolution strategy.
\newblock {\em Neural Computing and Applications}, 35:449--467, 2023.

\bibitem{Hecht2017}
Eugene Hecht.
\newblock {\em Optics}.
\newblock Pearson, Boston, 5th edition, 2017.

\end{thebibliography}
